\documentclass{ieja}
\usepackage{amsmath,amsthm,amssymb,amsfonts,amscd, array}
\def\currentvolume{}
\def\currentyear{}
\def\pagerange{}

\usepackage{pictex}

\DeclareMathOperator{\Ext}{Ext}
\DeclareMathOperator{\mo}{mod}
\DeclareMathOperator{\Hom}{Hom}

\def\mo{\operatorname{mod}}

\def\Hom{\operatorname{Hom}}

\def\Ext{\operatorname{Ext}}

\def\add{\operatorname{add}}

\newtheorem{thm}{Theorem}[section]

 \newtheorem{prop}[thm]{Proposition}

\title{The Ingalls-Thomas Bijections}
\author{M. A. A. Obaid, S. K. Nauman, W. M. Fakieh and C. M. Ringel}
\address{
Mustafa A. A. Obaid\\ 
E-mail: drmobaid@yahoo.com\\
\mbox{} \\
S. Khalid Nauman\\ 
E-mail: snauman@kau.edu.sa\\
\mbox{} \\
Wafaa  M. Fakieh\\ 
E-mail: wafaa.fakieh@hotmail.com\\
\mbox{} \\
Claus Michael Ringel\\ 
E-mail: ringel@math.uni-bielefeld.de\\
\mbox{} \\
King Abdulaziz University, Faculty of Science, \\
P.O.Box 80203, Jeddah 21589, Saudi Arabia
}

\begin{document}

\begin{abstract}
Given a finite acyclic quiver $Q$ with path algebra $\Lambda$,
Ingalls and Thomas have exhibited a bijection between the set of Morita equivalence classes of
support-tilting modules and the set of 
thick subcategories of $\mo\Lambda$ with covers, and they have collected
a large number of further bijections with these sets. We add some additional bijections
and show that all these bijections hold for arbitrary hereditary artin algebras.  
The proofs presented here seem to be of interest also in the special case of
the path algebra of a quiver. \\[+2mm]
 \subjclassname{\ Primary:
   16D90, 
   16G70. 
   Secondary:
   16G20,  
   05E10. 
  }
       Normal module, conormal module.  Torsion pair. Antichain. Exceptional antichain.
 \end{abstract}

\maketitle

\section{Introduction}

\subsection{}
Let $\Lambda$ be a  hereditary
artin algebra. We recall that an artin algebra $\Lambda$ is a $k$-algebra which is 
of finite length when considered as a $k$-module, where $k$ is a commutative artinian ring.
An artin algebra is {\it hereditary} provided submodules of projective modules are
projective. Since this means that  
the functors $\Ext_\Lambda^i$ vanish for $i \ge 2$, we write
$\Ext(M,M')$ instead of $\Ext_\Lambda^1(M,M').$ 
A typical example of a hereditary artin algebra is the path algebra of 
a finite acyclic quiver
(if $k$ is an algebraically closed field,  any 
hereditary artin $k$-algebra is Morita-equivalent to the path algebra of a finite
acyclic quiver, but otherwise there are many other hereditary artin $k$-algebras).

We will consider left $\Lambda$-modules of finite
length and call them just modules. The category of all modules will be denoted by $\mo\Lambda$.
We denote by $n = n(\Lambda)$ the {\it rank} of $\Lambda$; this is by definition
the number of isomorphism classes of 
simple modules. 

Given a module $M$, we denote by $\Lambda(M)$ its
{\it support algebra of} $M$; this is the factor algebra of $\Lambda$
modulo the ideal generated by all idempotents $e$ with $eM = 0.$ 
The support algebra $\Lambda(M)$ is again a hereditary
artin algebra (but usually not connected, even if $\Lambda$ is connected). The rank of the support algebra
of $M$ will be called the {\it support-rank} of $M$.
If $M$ is a module, the set of simple modules $S$ which occur as composition factors of $M$ 
will be called the {\it support} of $M.$ The module $M$ is said to be {\it sincere} provided
any simple module belongs to the support of $M$ (thus provided the only idempotent
$e \in \Lambda$ with $eM = 0$ is $e=0$). 

\subsection{}
The subcategories of $\mo\Lambda$ which we will consider are full subcategories
which are closed under direct sums and direct summands.
Given a class $\mathcal X$ of modules, we denote by 
$\add \mathcal X$ the class of modules which are direct summands of direct sums of modules in
$\mathcal X$. If $\mathcal X = \{X\}$ for a single module $X$, we write $\add X$ instead of
$\add\{X\}$.  
The modules $X, X'$ are said to be {\it Morita equivalent} provided $\add X = \add X'$.
Note that multiplicity-free modules which are Morita equivalent are actually isomorphic.
On the other hand, every module is Morita equivalent to a multiplicity-free module.
   
\subsection{Support-tilting modules.}Following earlier considerations of Brenner and Butler,
tilting modules have been defined in \cite{[HR]}. 
We say that a module $M$ has no self-extensions, provided $\Ext(M,M) = 0$. 
In the present setting, 
a module $T$ without self-extensions is said to be a {\it tilting module} provided it has 
precisely $n$ isomorphism classes of indecomposable direct summands (where $n$ is the
rank of the artin algebra $\Lambda$), or, equivalently,
provided ${}_\Lambda\Lambda$ is the kernel of a surjective map in $\add T$ (or an injective cogenerator
is the cokernel of an injective map in $\add T)$.
A module $M$ is said to be
{\it support-tilting} provided $M$ considered as a $\Lambda(M)$-module is a tilting module. 

Here is one of the sets we are interested in: {\it the set of Morita equivalence classes of
support-tilting modules.}

\subsection{Thick subcategories with a cover.}
A subcategory $\mathcal A$ of $\mo\Lambda$ is called a {\it thick} 
(or {\it wide}) subcategory provided it is
closed under kernels, cokernels and extensions. Note that a thick subcategory is an abelian category, and the
inclusion functor $\mathcal A \to \mo\Lambda$ is exact.

A module $X$ is said to 
{\it generate} a module $Y$ provided $Y$ is a factor module of a direct sum of copies of $X$.
Dually, a module $X$ {\it cogenerates} a module $Y$ provided $Y$ is a
submodule of a direct sum of copies of $X$ (since the modules considered here are of
finite length, it is sufficient to look at direct sums of copies of $X$; for general
modules one would have to use products). 
Given a class $\mathcal X$ of modules, let $\mathcal G(\mathcal X)$ be the subcategory
of all modules which are generated by modules in 
$\add\mathcal X$, and let $\mathcal H(\mathcal X)$ be the subcategory
of all modules which are cogenerated by modules in $\add\mathcal X$.
If $\mathcal C$ is a subcategory and $C\in \mathcal C$, then $C$ is said to be a {\it cover} of $\mathcal C$ provided
$\mathcal C \subseteq \mathcal G(C),$ and $C$ is said to be a {\it cocover} of $\mathcal C$ provided
$\mathcal C \subseteq \mathcal H(C).$ 

This is the second set of interest: {\it the set of
thick subcategories of $\mo\Lambda$ with covers.}
      
\subsection{}
If $\Lambda$ is the path algebra of a finite acyclic quiver, Ingalls and Thomas have exhibited
a bijection between the set of Morita equivalence classes of
support-tilting modules and the set of 
thick subcategories of $\mo\Lambda$ with covers.
The aim of this paper is to provide a proof of the Ingalls-Thomas bijection for 
arbitrary hereditary artin algebras.
Our proof draws attention to three additional sets which are in bijection with the set of
Morita equivalence classes of support-tilting modules:
the set of isomorphism classes of exceptional antichains in $\mo\Lambda$, 
as well as the set of isomorphism classes of normal or of conormal modules without self-extensions.
Here are the definitions.
     
\subsection{Exceptional antichains.} Given an additive category $\mathcal C$, a {\it brick} is $\mathcal C$ is
an object whose endomorphism ring is a division ring.  
Bricks $A_1, A_2$ are said to be 
{\it orthogonal,} provided $\Hom(A_1,A_2) = 0 = \Hom(A_2,A_1)$.
An {\it antichain} $A = \{A_1,\dots,A_t\}$ in $\mathcal C$
is a set of pairwise
orthogonal bricks  (antichains are called discrete subsets
in \cite{[GP]} and $\Hom$-free subsets in \cite{[HK]}, see also the remark \ref{antichain}). 
Antichains $A = \{A_1,\dots,A_t\}$ and $A' = \{A'_1,\dots,A'_{t'}\}$ are said to
be {\it isomorphic,} provided the objects $\bigoplus_i A_i$ 
and $\bigoplus_j A'_j$ are isomorphic.

Given an antichain $A = \{A_1,\dots,A_t\}$ in $\mo\Lambda$,
its $\Ext$-quiver $Q_A$ has as vertices the elements $A_i$
and there is an arrow $A_i\to A_j$ provided $\Ext(A_i,A_j) \neq 0$ (one may endow this quiver with
a valuation, taking into account the size of the $\Ext$-groups, but this is not needed here). 
We say that an antichain $A$ is {\it exceptional,} provided its $\Ext$-quiver $Q_A$ is acyclic, thus 
provided we may index the elements of $A$ in such a way that $\Ext(A_i,A_j) = 0$
for all pairs $i \ge j$.  

\subsection{Normal (or conormal) modules without self-extensions.}
A module $M$ is said to be {\it normal} provided 
given a direct decomposition $M = M'\oplus M''$ such that
$M'$ generates $M''$, we have $M'' = 0.$  And $M$ is {\it conormal} provided
given a direct decomposition $M = M'\oplus M''$ such that
$M'$ cogenerates $M''$, we have $M'' = 0.$

There is the following well-known fact (see, for example \cite{[R2]}):
{\it A sincere module without self-extensions is faithful,} thus any module $M$ without
self-extensions is a faithful $\Lambda(M)$-module.

\subsection{}
Since its introduction, tilting theory concerns the study of suitable torsion
pairs in $\mo\Lambda$. It seems worthwhile to include this aspect in our considerations. 
Recall that a {\it torsion class} in $\mo \Lambda$ is a class of 
modules which is closed under factor modules and extensions.
A {\it torsionfree class} in $\mo\Lambda$ is a class of modules which is 
closed under submodules and extensions.
       	     \bigskip
	     
It was the decisive idea of Ingalls and Thomas \cite{[IT]} to relate
the support-tilting modules to thick subcategories and to exhibit in this way
a number of bijections. They were dealing with path algebras of finite acyclic quivers,
here we consider the case of an arbitrary hereditary artin algebra.

\begin{thm}\label{theorem} Let $\Lambda$ be a hereditary artin algebra.
There are bijections between the following data:

\begin{itemize}
\item{\rm(1)} Isomorphism classes of exceptional antichains.
\item{\rm(2)} Thick subcategories with a cover.
\item{\rm(3)} Isomorphism classes of normal modules without self-extensions.
  \medskip
  
\item{\rm(4)} Morita equivalence classes of support-tilting modules.
\item{\rm(5)} Torsion classes with a cover.
\end{itemize}
	\medskip

	\noindent
  If $\Lambda$ is in addition representation-finite, then
  
\begin{itemize}
\item{\rm(1$'$)} {\it All antichains are exceptional.}
\item{\rm(2$'$)} {\it All thick subcategories have a cover.}
\item{\rm(5$'$)} {\it All torsion classes have a cover.}
\end{itemize}
\end{thm}
		 
We have separated the five sets in Theorem \ref{theorem} into two groups, since 
there is a great affinity between (1), (2) and (3) on the one hand, and (4) and (5)
on the other hand. 
The essential bijection concerns the sets (2) and (4). As we have mentioned,
such a bijection was exhibited  by Ingalls-Thomas \cite{[IT]}
in case $\Lambda$ is the path algebra of a finite acyclic quiver. 
A bijection
between (4) and (5) has been known for a long time.
A bijection between (1) and (2) was
exhibited already in 1976, see \cite{[R1]}. For a bijection between (1) and (3), one may 
refer to \cite{[DR]}, as we will see below.
      
\subsection{Outline of the paper.}
Sections 2 to 4 provide the required bijections in detail,
and an outline of the corresponding proofs. 
Section 5 is devoted to duality. Whereas the sets of the
form (1), (2) and (4) are preserved under duality, this is not the case for the sets (3) and (5),
thus, using duality, we obtain bijections with two further sets: the set (6) 
of {\it isomorphism classes of conormal modules without self-extensions,} and the set
(7) of the {\it torsionfree classes with a cocover.}
In the final section 6 we deal with the support of the various modules and subcategories.

As a supplement of the theorem, we 
have mentioned that for $\Lambda$ representation-finite, certain conditions
are always satisfied. First of all, 
if $\Lambda$ is representation-finite, then any subcategory of $\mo\Lambda$ has
both a cover and a cocover. And second, 
it is well-known that for an antichain $A$ which is not exceptional, the class $\mathcal F(A)$ 
of all modules with a filtration with factors in $A$
contains infinitely many isomorphism classes of
indecomposable $\Lambda$-modules, thus $\Lambda$ cannot be representation-finite.
	       
\subsection{}
The case of $\Lambda$ being representation-finite is studied in more detail
in our paper \cite{[ONFR]}. Such an artin algebra $\Lambda$
is called a {\it Dynkin algebra}, since the underlying graph 
of its valuated quiver is the disjoint union of Dynkin diagrams. 
There, we will discuss the number
of tilting and support-tilting modules for these algebras.
For the Dynkin cases $\mathbb A$, we obtain the Catalan triangle,
for the cases $\mathbb B$ and $\mathbb C$ we obtain the increasing part of the
Pascal triangle, and finally for the cases $\mathbb D$ we obtain an 
expansion of the increasing part of the Lucas triangle.    
For a further study of the Ingalls-Thomas bijections in general, we also 
may refer to the forthcoming survey \cite{[R4]}. 
    
\section{The bijections between (1), (2) and (3)}

From (1) to (2): {\it If $A$ is an antichain, take $\mathcal F(A)$,} this is the set of all
$\Lambda$-modules with a filtration with factors in $A.$ 
The full subcategory $\mathcal F(A)$ is an abelian category 
with exact embedding functor
and obviously closed under extensions, 
its simple objects are just the elements of $A$; the process of
considering the elements of 
$A$ as objects in $\mathcal F(A)$ is called {\it simplification} in \cite{[R1]}. 
If the antichain $A$ is
exceptional, the category $\mathcal F(A)$ is equivalent to the module category of an artin
algebra, thus it has projective generators. 
Every projective generator of $\mathcal F(A)$ is a cover for $\mathcal F(A)$.

For the step (1) to (2), we also may refer to \cite{[DR]}. Namely,
an exceptional antichain $A$ is a standardizable set as considered in \cite{[DR]} and the proof
of Theorem 2 in \cite{[DR]}
asserts that there is a quasi-hereditary algebra $B$ such that the
subcategory $\mathcal F(A)$ is equivalent to the category of $\Delta$-filtered $B$-modules. Since the
standardizable set $A$ consists of pairwise orthogonal modules, the same is true 
for the $\Delta$-modules
of $B$, and consequently the $\Delta$-modules of $B$ are just the simple $B$-modules. This shows that
the category of $\Delta$-filtered $B$-modules is the whole category $\mo B.$

From (2) to (1): {\it If $\mathcal A$ is a thick subcategory with a cover, let $\mathcal S(\mathcal A)$ 
be the set of simple
objects in $\mathcal A$,} one from each isomorphism class. Then 
$\mathcal S(\mathcal A)$ is an exceptional antichain in $\mo\Lambda$. 
Namely, a thick subcategory with a cover is equivalent, as a category, to the module
category $\mo \Lambda'$ of an artin algebra $\Lambda'$. Such an equivalence identifies 
the quiver $Q_{\mathcal S(\mathcal A)}$
with the quiver of the artin algebra 
$\Lambda'$ (the quiver of an artin algebra is just the $\Ext$-quiver
of the simple $\Lambda'$-modules). It is well-known (and easy to see) that the quiver
of an artin algebra is acyclic. 

From (2) to (3): {\it If $\mathcal A$ is a thick subcategory with a cover, let $P$ be a minimal projective generator of
$\mathcal A$.} Then $P$ is a normal module without self-extensions.

If we start with (1), say with an exceptional
antichain $A$, and use \cite{[DR]} in order to find an equivalence $\eta\!:\mathcal F(A) \to \mo B$,
the proof of
Theorem 2 in \cite{[DR]} first constructs indecomposable objects in $\mathcal F(A)$ which correspond under $\eta$ to the
indecomposable projective $B$-modules. In this way, one constructs a minimal projective generator for the
abelian category $\mathcal F(A)$.

From (3) to (1). Let $N$ be a normal module without self-extensions. Write $N = \bigoplus_i N_i$ with indecomposable
modules $N_i$. For any $i$, let $u_i\!:U_i \to N_i$ be a minimal right $\mathcal N_i$-approximation of $N_i$, where
$\mathcal N_i = \add(\{N_j\mid j\neq i\}$. Since $N$ is normal, the map $u_i$ cannot be surjective. Since $\Lambda$ is
hereditary, it follows that $u_i$ is injective
and we denote by $p_i\!:N_i \to \Delta(i)$ the cokernel of $u_i$.
Since $u_i$ is not surjective, we see that $\Delta(i) \neq 0.$
We claim that the modules $\Delta(i)$ are pairwise orthogonal bricks. Let $h\!:N_j \to \Delta(i)$ be a map, and form the
induced exact sequence
$$
{\beginpicture
\setcoordinatesystem units <2cm,1.3cm>
\multiput{$0$} at 0 0  4 0  0 -1  4 -1 /  
\put{$U_i$} at 1 0  
\put{$M$} at 2 0  
\put{$N_j$} at 3 0 
\put{$U_i$} at 1 -1 
\put{$N_i$} at 2 -1  
\put{$\Delta(i)$} at 3 -1  
\arrow <1.5mm> [0.25,0.75] from 0.3 0 to 0.7 0
\arrow <1.5mm> [0.25,0.75] from 1.3 0 to 1.7 0
\arrow <1.5mm> [0.25,0.75] from 2.3 0 to 2.7 0
\arrow <1.5mm> [0.25,0.75] from 3.3 0 to 3.7 0
\arrow <1.5mm> [0.25,0.75] from 0.3 -1 to 0.7 -1
\arrow <1.5mm> [0.25,0.75] from 1.3 -1 to 1.7 -1
\arrow <1.5mm> [0.25,0.75] from 2.3 -1 to 2.7 -1
\arrow <1.5mm> [0.25,0.75] from 3.3 -1 to 3.7 -1
\plot 1 -.3  1 -.7 /
\plot 1.05 -.3  1.05 -.7 /
\arrow <1.5mm> [0.25,0.75] from 2 -.3 to 2 -.7
\arrow <1.5mm> [0.25,0.75] from 3 -.3 to 3 -.7
\put{$h$} at 3.13 -.45
\endpicture}
$$
Since $U_i$ belongs to $\mathcal N_i$ and $N$ has no self-extensions, we have $\Ext(N_j,U_i) = 0$,
thus the upper sequence splits. It follows that  
there is a map $h'\!:N_j \to N_i$ such that $h = p_ih'.$ 
This has two consequences.

First of all, consider the case $j = i$. Let $g$ be any
endomorphism of $\Delta(i)$ and look at the map $h = gp_i\!:N_i \to \Delta(i)$.
We see that there is an endomorphism $g'\!:N_i \to N_i$
with $gp_i = p_ig'$. Since all non-zero endomorphisms of $N_i$ are invertible, 
the same is true for $\Delta(i)$.
In this way, we see that $\Delta(i)$ is a brick. 

Second, let $g\!:\Delta(j) \to \Delta(i)$ be a homomorphism
with $j\neq i$ and consider $h = gp_j\!:N_j \to \Delta(i).$ 
There is $g'\!:N_j \to N_i$ such that $gp_j = p_ig'$. Since $u_i$ is a left
$\mathcal N_i$-approximation, it follows that $g' = u_ig''$ for some $g''\!:N_j \to U_i.$ But then $gp_j = p_ig' =
p_iu_ig'' = 0$ and therefore $g = 0.$

In this way, we have shown that $\Delta = \{\Delta(i)\mid i\}$ is an antichain.
Using induction on the length $|N_i|$ 
of $N_i$, we see that $N_i$ belongs to $\mathcal F(\Delta).$ Namely, if $N_i$ is of length
$1$, then $U_i = 0$ since $\Delta(i) \neq 0.$ If $|N_i| \ge 2,$ then $U_i$ is a direct sum of modules of the form
$N_j$ with $|N_j| < |N_i|$, thus by induction $U_i$ belongs to $\mathcal F(\Delta)$ and therefore also $N_i$ belongs to
$\mathcal F(\Delta)$.

The surjective map $p_i\!:N_i\to \Delta(i)$ yields a surjective map $\Ext(N,N_i) \to \Ext(N,\Delta(i))$, thus
$\Ext(N,\Delta(i)) = 0$ for all $i$, and therefore $\Ext(N,M) = 0$ for all $M \in \mathcal F(\Delta)$.
This shows that the objects $N_i$ are indecomposable projective objects in $\mathcal F(\Delta)$; actually, $N_i$
is the projective cover of $\Delta(i)$ in $\mathcal F(\Delta).$ As usual, one sees now that $\Ext(\Delta(i),\Delta(j)) \neq 0$
if and only if $N_j$ is a direct summand of $U_i$. If $N_j$ is a direct summand of $U_i$, then, in particular,
$|N_j| < |N_i|$. This shows that the $\Ext$-quiver of $\Delta$ is acyclic.

Starting with an exceptional antichain $A$ in (1), and going via (2) to (3), we obtain a 
minimal projective
generator $P$ of $\mathcal F(A)$. Going from (3) to (1), we attach to $P$ the antichain $\Delta$ whose elements
are just the simple objects in $\mathcal F(A)$, but these are just the elements of $A$.
Conversely, starting in (3) say with a normal module $N$ without self-extensions, then going to (1), we attach to it the
antichain $\Delta$. Going via (2) to (3), we form a minimal projective generator in $\mathcal F(\Delta).$ But $N$
is up to isomorphism the only minimal projective generator in $\mathcal F(\Delta).$

\section{The bijection between (3) and (4)}\label{3}

From (4) to (3): {\it If $T$ is a support-tilting module, let $\nu(T)$ be its normalization.} This clearly is
a normal module without self-extensions. Here we use that any module $M$ can be written in the form $M = M'\oplus M''$
where $M'$ is normal and generates $M''$ (this of course is trivial), and that such a decomposition is
unique up to isomorphism (this is not so obvious); the module $M'$ is called a {\it normalization} of the
module $M$. The uniqueness was first shown by Roiter \cite{[Ro]} 
and then also by Auslander-Smal\o{} \cite{[AS]}, see also
\cite{[R3]}. The uniqueness shows that the map $\nu$ going from (4) to (3) is well-defined.

Let us show that $\nu$ is injective when we are dealing with
support-tilting modules. We claim the following:
if $T, T'$ are support-tilting modules with $\nu(T) =
\nu(T'),$ then $T$ and $T'$ are Morita equivalent. For the proof, we may replace $\Lambda$
by the support algebra $\Lambda(T) = \Lambda(T')$, thus we 
may assume that $T, T'$ are tilting modules. Now,
$T'$ is generated by $\nu(T') = \nu(T)$, thus
by $T$. Since $T$  generates $T'$, it follows from $\Ext(T,T)=0$ that $\Ext(T,T') = 0.$
Similarly,  $T'$ generates $T$ and therefore $\Ext(T',T) = 0.$ Altogether we see that
$\Ext(T\oplus T',T\oplus T') = 0.$ Since $T$ is a tilting module, this implies that $T'$ belongs to $\add T$.
Similarly, since $T'$ is a tilting module, we see that $T$ belongs to $\add T'$.

In order to see that $\nu$ is
also surjective, we need to find for any normal module $N$ without self-extensions a support-tilting module $T$
with $\nu(T) = N.$ This we will show next.

From (3) to (4): If $N$ is  a module without self-extensions,
there is a module $Y$, with the
following properties: first, $Y$ is generated by $N$, and 
second, $N\oplus Y$ is a support-tilting module; we call $Y$ 
a {\it factor complement} for  $N$
(this is the dual version of forming a Bongartz complement, see for example \cite{[R2]}). 

Here is the construction of a factor complement $Y$ of a module without self-extensions
(we follow \cite{[R2]}).
Let $\Lambda(N)$ be the support algebra for $N$ and $Z$ an injective cogenerator
for $\mo \Lambda(N)$. 
We claim that {\it there exists an epimorphism 
$Y \to Z$ with kernel in $\add N$ such that $\Ext(Y,N) = 0$.} Such an epimorphism
can be obtained as a universal foundation of $Z$ by $N$
(sometimes also called a universal extension of $Z$ by $N$ from below): take
exact sequences $0 \to N \to Y_i \to Z \to 0$ such that the corresponding elements in $\Ext(Z,N)$
generate it as a $k$-module, 
and form the direct sum of these sequences. The induced sequence with respect to
the diagonal inclusion $u\!:Z \to \bigoplus_i Z$
$$
{\beginpicture
\setcoordinatesystem units <2cm,1.3cm>
\multiput{$0$} at 0 0  4 0  0 -1  4 -1 /  
\put{$\bigoplus_i N$} at 1 0  
\put{$\bigoplus_i Y_i$} at 2 0  
\put{$\bigoplus_i Z$} at 3 0 
\put{$\bigoplus_i N$} at 1 -1 
\put{$Y$} at 2 -1  
\put{$Z$} at 3 -1  
\arrow <1.5mm> [0.25,0.75] from 0.3 0 to 0.7 0
\arrow <1.5mm> [0.25,0.75] from 1.3 0 to 1.7 0
\arrow <1.5mm> [0.25,0.75] from 2.3 0 to 2.7 0
\arrow <1.5mm> [0.25,0.75] from 3.3 0 to 3.7 0
\arrow <1.5mm> [0.25,0.75] from 0.3 -1 to 0.7 -1
\arrow <1.5mm> [0.25,0.75] from 1.3 -1 to 1.7 -1
\arrow <1.5mm> [0.25,0.75] from 2.3 -1 to 2.7 -1
\arrow <1.5mm> [0.25,0.75] from 3.3 -1 to 3.7 -1
\plot 1 -.3  1 -.7 /
\plot 1.05 -.3  1.05 -.7 /
\arrow <1.5mm> [0.25,0.75] from 2 -.7 to 2 -.3
\arrow <1.5mm> [0.25,0.75] from 3 -.7 to 3 -.3
\put{$u$} at 3.13 -.55
\endpicture}
$$
yields a universal foundation $g\!:Y \to Z$. 

In general, {\it given a universal foundation
$g\!:Y   \to Z$ of $Z$ by $N$, say with kernel $N'$, the module $Y$ is generated by $N$.} 
Namely, since $N$ has no self-extensions, it is a
faithful $\Lambda(N)$-module, thus $Z$ is generated by $N$. An epimorphism
$h\!:N^t \to Z$ yields a commutative diagram with exact rows
$$
{\beginpicture
\setcoordinatesystem units <2cm,1.3cm>
\multiput{$0$} at 0 0  4 0  0 -1  4 -1 /  
\put{$N'$} at 1 0  
\put{$Y$} at 2 0  
\put{$Z$} at 3 0 
\put{$N'$} at 1 -1 
\put{$N''$} at 2 -1  
\put{$N^t$} at 3 -1  
\arrow <1.5mm> [0.25,0.75] from 0.3 0 to 0.7 0
\arrow <1.5mm> [0.25,0.75] from 1.3 0 to 1.7 0
\arrow <1.5mm> [0.25,0.75] from 2.3 0 to 2.7 0
\arrow <1.5mm> [0.25,0.75] from 3.3 0 to 3.7 0
\arrow <1.5mm> [0.25,0.75] from 0.3 -1 to 0.7 -1
\arrow <1.5mm> [0.25,0.75] from 1.3 -1 to 1.7 -1
\arrow <1.5mm> [0.25,0.75] from 2.3 -1 to 2.7 -1
\arrow <1.5mm> [0.25,0.75] from 3.3 -1 to 3.7 -1
\plot 1 -.3  1 -.7 /
\plot 1.05 -.3  1.05 -.7 /
\arrow <1.5mm> [0.25,0.75] from 2 -.7 to 2 -.3
\arrow <1.5mm> [0.25,0.75] from 3 -.7 to 3 -.3
\put{$h'$} at 2.13 -.55
\put{$h$} at 3.13 -.55
\endpicture}
$$
Since $\Ext(N,N) = 0$, the lower sequence splits, thus $N''$ belongs to $\add N$.
Since $h$ is surjective, also $h'$ is surjective, thus $Y$ is generated by $N$.

It remains to be seen that $N\oplus Y$ is support-tilting.
Since $N$ generates $Y$, it follows from $\Ext(N\oplus Y,N) = 0$ that 
$\Ext(N\oplus Y, Y) = 0.$ In this way, we see that $N\oplus Y$ has no self-extensions.
The exact sequence $0 \to \bigoplus_i N \to Y \to Z \to 0$ shows that $Z$ is the cokernel
of an injective map in $\add (N\oplus Y)$, thus $N\oplus Y$ is a support-tilting module.
This completes the proof that $Y$ is a factor complement for $N$.

If we choose a minimal direct summand 
$\phi(N)$ of $Y$ such that $N\oplus \phi(N)$ is a support-tilting module, then
$\phi(N)$ is uniquely determined by $N$ and may be called a {\it minimal
factor complement} for $N$. Thus, going from (3) to (4), we may attach to 
a normal module $N$ without self-extension the multiplicity-free support-tilting module
$N\oplus \phi(N).$

Of course, if $N$ is normal, then $N$ is the normalization of $N\oplus Y$.
Thus starting with a normal module $N$ without self-extensions, then
going from (3) to (4) and back to (3), we obtain $N$. On the other hand, let $T$ be
support-tilting. From (4) to (3) we take $\nu(T)$. From (3) to (4), we add to $\nu(T)$
a factor complement, say $N'$. But $T$ and $T' =\nu(T)\oplus N'$ both are support-tilting
modules with $\nu(T) = \nu(T')$ and generated by this module $\nu(T)$, 
thus they are Morita equivalent.

\section{The bijection between (4) and (5)}

First, we show the following: If $T$ is a support-tilting module and $\mathcal G = \mathcal G(T)$, then
$\add T$ is the class of the $\Ext$-projective modules in $\mathcal G$.
Tilting theory asserts that $\mathcal G$ is the class of $\Lambda(T)$-modules $M$ such that $\Ext(T,M) = 0$.
Let $M$ be in $\mathcal G$ and $g\!:T' \to M$ be a right $T$-approximation of $M$. Then $g$ is
surjective and the kernel $M'$ of $g$ satisfies $\Ext(T,M') = 0$, thus belongs to $\mathcal G$.
If $M$ is $\Ext$-projective, then the exact sequence $0 \to M' \to T' \to M \to 0$ splits, thus $M$ is in
$\add T.$ This shows that the $\Ext$-projective modules in $\mathcal G$ are just the modules in $\add T.$

From (4) to (5): If $T$ is a module without self-extensions, let $\mathcal G(T)$ be the class of modules
generated by $T$.
Then it is well-known (and easy to see) 
that $T$ is a torsion class. Of course, $T$ is a cover for $\mathcal G(T)$.

From (5) to (4): If $\mathcal C$ is a torsion class with a cover $C$, then we attach to it a module $T$
such that $\add T$ is the class of $\Ext$-projective modules in $\mathcal G$. In order to do so, we need to know
that the class $\mathcal E$ of $\Ext$-projective modules
in $\mathcal C$ is finite, say $\mathcal E = \add T$ for some module $T$. We also have to show that $T$
is support-tilting.

Along with $C$, 
its normalization $\nu(C)$ is also a cover. A normal cover of a torsion class has no self-extension
(see Proposition 1 of \cite{[R3]}). Let $B$ be a factor complement for $\nu(C).$ As we have seen,
$T = \nu(C)\oplus B$ is a support-tilting module. Since $B$ is generated by $\nu(C)$, we have $\mathcal G(T) =
\mathcal G(\nu(C)) = \mathcal G(C) = \mathcal C.$ But we have shown already that $\add T$ is the class of
$\Ext$-projective modules in $\mathcal G(T)$.

From (4) to (5) to (4): Let us start with a support-tilting module $T$ and attach to it $\mathcal G = \mathcal G(T)$.
As we have seen, the class of $\Ext$-projectives in $\mathcal G$ is $\add T$. We choose $T'$ with $\add T' = \add T$.
But this just means that $T, T'$ are Morita equivalent.

From (5) to (4) to (5). We start with a torsion class $\mathcal C$ with a cover, we choose a support-tilting module $T$
with $\mathcal C = \mathcal G(T),$  thus we are back at $\mathcal C$.
      
\section{Duality}

By definition, given an artin algebra $\Lambda$, there is a commutative artinian ring $k$
such that $\Lambda$ is a $k$-algebra and is of finite length when considered as a
$k$-module.  If $\Lambda$ is an artin algebra, also the opposite algebra $\Lambda^{\text{op}}$
is an artin algebra. If we denote by $E$ a minimal injective cogenerator for $\mo k$,
the functor $D = \Hom_k(-,E)$ provides an equivalence between $\mo\Lambda$ and
$(\mo\Lambda^{\text{op}})^{\text{op}}.$ We can use this duality in order to exhibit further
bijections. 

Using duality, the sets (1), (2) and (4) are preserved. 
Of course, the dual concept of a thick subcategory with
a cover is a thick subcategory with a cocover.
An abelian $k$-category with finitely many simple objects and such that the
$\Hom$ and $\Ext$-groups are $k$-modules of finite length, 
has a cover if and only if it has a cocover.
    \smallskip

    Dualizing (3) we get:
    	      
\begin{itemize}
\item{(6)} {\it The isomorphism classes of conormal modules without self-extensions.}
\end{itemize}

Dualizing (5) we get:

\begin{itemize}
\item{(7)} {\it The torsionfree classes with a cocover.}
\end{itemize}
	\medskip 

{\it The sets defined in {\rm(6)} and {\rm(7)} correspond bijectively to the sets}
(1),\dots,(5).
       \bigskip 

{\bf Remark.} The bijections between the set (2) of thick subcategories $\mathcal A$ 
and the sets (1), (3) and (6) of isomorphism classes of suitable modules 
can be reformulated as follows:
In an abelian category we may look at the semi-simple, the
projective and the injective objects: the set of simple objects in $\mathcal A$ is
an antichain in $\mo\Lambda$, a minimal projective generator in $\mathcal A$
is a normal module without self-extensions, a minimal injective cogenerator
is a conormal module without self-extensions.  
These are the procedures to obtain from a thick subcategory the corresponding antichain,
as well as a normal or conormal module without self-extensions.

Conversely, let us start with (1), (3) or (6). It has been
mentioned already that starting with an antichain $A$, we take the full subcategory $\mathcal F(A)$
of all modules with a filtration with factors in $A$. Starting with a normal module $P$
without self-extensions, the corresponding thick subcategory $\mathcal A$ consists of all modules
which arise as the cokernel of a map in $\add P$ (in this way, we specify projective
presentations of the objects in $\mathcal A).$ Dually, starting with a conormal module $I$
without self-extensions, the corresponding thick subcategory $\mathcal A$ consists of all modules
which arise as the kernel of a map in $\add I$ (in this way, we specify injective
presentations of the objects in $\mathcal A).$ 
	      
\section{The support of a module, sincere modules and subcategories}

\begin{prop} The bijections which we have constructed preserve the support.
\end{prop}

Specializing the Ingalls-Thomas bijections to sincere modules, it follows from 
the proposition that we get
{\it bijections between:}
    
\begin{itemize}
\item{(1)} {\it Isomorphism classes of exceptional sincere antichains.}
\item{(2)} {\it Thick subcategories with a sincere generator.}
\item{(3)} {\it Isomorphism classes of normal sincere modules without self-extensions.}
\item{(4)} {\it Morita equivalence classes of tilting modules.}
\item{(5)} {\it Torsion classes with a sincere generator.}
\item{(6)} {\it Isomorphism classes of conormal sincere modules without self-extensions.}
\item{(7)} {\it Torsionfree classes with a sincere cogenerator.}
\end{itemize}

Of course, conversely this special case implies the general case.

\section{Final remarks}

\subsection{}  The aim of our discussion was to extend results of Ingalls and Thomas
which were established for path algebras of finite acyclic quivers to arbitrary 
hereditary artin algebras. 
Experts may not be surprised that results concerning path algebras of finite
acyclic quivers can be extended in this way: after all,  there is a general feeling that such 
generalizations are always possible. But the paper \cite{[ORT]} may serve as a warning. The paper 
provides a description of the cofinite quotient-closed subcategories of $\mo\Lambda$,
where $\Lambda$ is the path algebra of a finite acyclic quiver. In section 9 of \cite{[ORT]},
the author discuss the problem of extending the result to finite-dimensional hereditary 
$k$-algebras, but they are able to provide a solution only in the case of $k$ being a finite field. 

On the other hand, one may ask whether the setting may be further
enlarged to deal with hereditary artinian or even hereditary semi-primary rings, and not just with
hereditary artin algebras. 
Note that our considerations use duality arguments and finiteness conditions which
rely on the artin algebra assumption. 

\subsection{}  A further possible generalization has been stressed by the referee: to drop the condition on
$\Lambda$ to be hereditary, thus to deal with an arbitrary artin algebra. 
For any finite-dimensional $k$-algebra $\Lambda$, with $k$ an algebraically closed field, the paper 
\cite{[AIR]} by Adachi, Iyama and Reiten provides a bijection between support $\tau$-tilting modules
in $\mo \Lambda$ and torsion classes with covers, extending in this way 
the corresponding Ingalls-Thomas
bijection (for a hereditary artin algebra, the $\tau$-tilting modules are just the
tilting modules). 
Also, let us remark that the relationship between torsion 
classes and thick
subcategories in $\mo\Lambda$ has been discussed by Marks and Stovicek \cite{[MS]}.

\subsection{}\label{antichain}  Our presentation of the Ingalls-Thomas bijections is centered around the notion
of antichains in additive categories. 
Let us motivate the definition. Given a poset $P$, a {\it chain} in $P$ is a 
subset of pairwise comparable elements, 
whereas an {\it antichain} in $P$ is a subset of pairwise incomparable elements. 
Now consider the linearization $kP$ of $P$, were $k$ is a field: 
this is an additive $k$-category whose indecomposable
objects are the elements of $P$ such that $\Hom_{kP}(x,y) = k$ 
provided $x \le y$ in $P$ and $\Hom_{kP}(x,y) = 0$
otherwise, such that the composition of maps in  $kP$ 
is given  by the multiplication in $k$, and, finally, such that any object in
$kP$ is the direct sum of indecomposable objects. Of course, a subset
$A$ of $P$ is an antichain in $P$ if and only if  $A$ (considered as a set of objects in
$kP$) consists of pairwise orthogonal bricks 
(thus, is an antichain in the additive category $kP$). As we see, antichains in additive categories
have to be considered as a direct generalization of antichains in posets.

The reader should be aware that starting with a Dynkin diagram $\Delta$
and its set $\Phi_+(\Delta)$ of positive roots, 
several kinds of (different, but related) antichains have to be distinguished: 
First of all,
$\Phi_+(\Delta)$ is in an intrinsic way a poset, called the {\it root poset} of type $\Delta$,
and we
may consider the set $\bold A(\Delta)$ of antichains in this root poset $\Phi_+(\Delta)$.
Second, choosing an orientation $\Omega$ of the
Dynkin diagram (or, equivalently, a Coxeter element in the corresponding
Weyl group), we may identify the elements of $\Phi_+(\Delta)$ with the 
indecomposable $\Lambda$-modules, thus with the indecomposable 
objects in the additive category $\mo\Lambda$. 
The set  
of antichains in $\mo\Lambda$ only depends on $\Delta$ and $\Omega$
(and not on the choice of $\Lambda$), thus we may denote it by $\bold A(\Delta,\Omega)$.
It is known for a long time that 
the set $\bold A(\Delta)$ of antichains in the root poset $\Phi_+(\Delta)$ 
and the set $\bold A(\Delta,\Omega)$ of 
antichains in $\mo\Lambda$ have the same enumeration
(for a uniform proof, see \cite{[AST]}), 
but a fully satisfactory explanation is still missing.
In the case of the quiver $\Bbb A_n$ with linear orientation, 
this concerns the quite obvious
bijection between non-nesting and non-crossing partitions. Note that
if $\Omega$ and $\Omega'$ are orientations of $\Delta$, it is easy to construct a
natural bijection between $\bold A(\Delta,\Omega)$ and $\bold A(\Delta,\Omega')$.
For a detailed discussion of the sets $\bold A(\Delta)$ and $\bold A(\Delta,\Omega)$,
we may refer to \cite{[R4]}. 

   \bigskip

\noindent {\bf Acknowledgment.} The authors want to thank the referee for very helpful
comments, in particular for spotting a wrong argument in section \ref{3}.
This work is funded by the Deanship of Scientific Research,
King Abdulaziz University, under grant No. 2-130/1434/HiCi.
The authors, therefore, acknowledge technical and financial support of KAU.

\end{document}